\documentclass[11pt]{article}

\makeindex
\usepackage[all]{xy}
\usepackage{graphicx}
\usepackage{amsmath}
\usepackage{amsfonts}
\usepackage{amssymb}
\usepackage{amsthm}

\newcommand{\Z}{{\bf{Z}}}
\newcommand{\Q}{{\bf{Q}}}
\newcommand{\Qbar}{{\overline{\Q}}}

\newcommand{\C}{{\bf{C}}}
\newcommand{\F}{{\bf{F}}}
\newcommand{\T}{{\bf{T}}}

\newcommand{\m}{{\mathfrak{m}}}

\newcommand{\tensor}{\otimes}

\newcommand{\ra}{\rightarrow}

\newcommand{\hra}{\hookrightarrow}

 1

\newcommand{\phie}{\phi_{\scriptscriptstyle E}}

\newcommand{\Fp}{{{\F}_p}}

\newcommand{\nel}{{n_{\scriptscriptstyle{E}}}}
\newcommand{\nele}{{\tilde{n}_{\scriptscriptstyle{E}}}}
\newcommand{\rA}{{r_{\scriptscriptstyle{A}}}}
\newcommand{\rAs}{{r_{\scriptscriptstyle{A}}^2}}
\newcommand{\rAfs}{{r_{\scriptscriptstyle{A_f}}^2}}

\newcommand{\rAe}{{\tilde{r}_{\scriptscriptstyle{A}}}}
\newcommand{\rAfe}{{\tilde{r}_{\scriptscriptstyle{A_f}}}}
\newcommand{\rAf}{{r_{\scriptscriptstyle{A_f}}}}

\newcommand{\nA}{{n_{\scriptscriptstyle{A}}}}
\newcommand{\nAe}{{\tilde{n}_{\scriptscriptstyle{A}}}}
\newcommand{\nAfe}{{\tilde{n}_{\scriptscriptstyle{A_f}}}}
\newcommand{\nAf}{n_{\scriptscriptstyle{A_f}}}

\newcommand{\Afdual}{{A^{\vee}_f}}

\newcommand{\Hom}{{\rm Hom}}
\newcommand{\OO}{{\mathcal{O}}}

\newcommand{\End}{{\rm End}}

\newcommand{\isom}{\cong}

\newcommand{\Frobl}{{{\rm Frob}_\ell}}

\newcommand{\Zlp}{{\mathbf{Z}_{p}}}   

\usepackage{bm}
\bmdefine\bmu{\mu}

\newcommand{\ndiv}{{\not|\ }}

\newcommand{\Ihat}{{\widehat{I}}}
\newcommand{\IAhat}{{\widehat{I_A}}}

\newcommand{\comment}[1]{}

\newtheorem{lem}{Lemma}[section]
\newtheorem{cor}[lem]{Corollary}
\newtheorem{prop}[lem]{Proposition}

\newtheorem{thm}[lem]{Theorem}

\theoremstyle{definition}

\newtheorem{eg}[lem]{Example}

\DeclareMathOperator{\Ann}{Ann}

\DeclareMathOperator{\ord}{ord}
\DeclareMathOperator{\Jac}{Jac}

\newcommand{\thetitle}
{The Modular number, the Congruence number, and Multiplicity One}

\begin{document}

\parindent=2em

\title{\thetitle}
\author{Amod Agashe\footnote{The 
author was partially supported by National Science Foundation Grant No. 0603668. }
}
\date{}
\maketitle

\noindent Department of Mathematics, Florida State University, U.S.A.

\noindent {\em Email:} agashe@math.fsu.edu

\begin{abstract}
Let $N$ be a positive integer and let $f$ be a newform of weight~$2$
on~$\Gamma_0(N)$.
In a joint paper with K.~Ribet and W.~Stein, the authors introduced
the notions of 
the modular number 
and the congruence number of the  quotient
abelian variety~$A_f$ of~$J_0(N)$ associated
to the newform~$f$. These invariants are
analogs of the notions of the modular degree
and congruence primes respectively associated to optimal elliptic curves. 
In this article, we show that if $p$ is a prime number such that every maximal ideal
of the Hecke algebra of residue characteristic~$p$
that contains the annihilator ideal of~$f$
satisfies multiplicity one, 
then the modular number of~$A_f$ and the congruence number of~$A_f$ 
have the same $p$-adic valuation. We also discuss a more general setting
that involves certain abelian subvarieties of~$J_0(N)$ or~$J_1(N)$ and 
state a result about the structure of the intersection of
such abelian subvarieties as a module over the Hecke algebra, from which
the statement in the previous sentence follows.
We also give a numerical example where
our result implies the failure of
multiplicity one. 
\end{abstract}



\section{Introduction and some of the results}
\label{sec:intro}

Let $N$ be a positive integer and let $X_0(N)$ denote the modular curve
over~$\Q$ associated to the classification
of  isomorphism classes of elliptic curves with
a cyclic subgroup of order~$N$.  The Hecke algebra~$\T$ of level~$N$
is the subring of the ring of endomorphisms of $J_0(N)=\Jac(X_0(N))$
generated by the Hecke operators $T_n$ for all $n\geq 1$.  Let~$f$ be
a newform of weight~$2$ for~$\Gamma_0(N)$  and let $I_f$ denote 
${\rm Ann}_\T(f)$.
Then the quotient abelian variety
$A_f = J_0(N)/I_f J_0(N)$ is called the
newform quotient associated to~$f$.
If $f$ has integer Fourier coefficients, then $A_f$ is an elliptic
curve and by~\cite{breuil-conrad-diamond-taylor} 
any elliptic curve over~$\Q$
is isogenous to such an elliptic curve for some~$f$. 
The dual abelian variety~$\Afdual$ of~$A_f$ may be viewed
as an abelian subvariety of~$J_0(N)$. 
Recall that the {\em exponent} of a finite group~$G$ is the smallest positive
integer~$n$ such that multiplication by~$n$ annihilates
every element of~$G$.

The exponent of the group $\Afdual \cap I_f J$ is called
the {\em modular exponent} of~$A_f$, denoted~$\nAfe$, and its order is called
the {\em modular number}, denoted~$\nAf$
(see~\cite[\S3]{ars:moddeg}). Suppose for the moment
that $f$ has integer Fourier coefficients,
so that 
$A_f$ is an elliptic curve, which we 
denote by~$E$ for emphasis.
Composing the embedding $X_0(N)\hra J_0(N)$ that sends $\infty$ to~$0$
with the quotient map $J_0(N) \ra E$, we obtain a surjective morphism
of curves $\phie: X_0(N) \ra E$, whose degree is called
the {\em modular degree} of~$E$.
The modular exponent~$\nele$ of~$E$ is equal to its modular degree,
and the modular number~$\nel$ is the square of the modular degree
(see~\cite[\S3]{ars:moddeg}). 
In general, for any newform~$f$, the modular number~$\nAf$ is a perfect
square (e.g., see~\cite[Lemma~3.14]{agst:bsd}).

Let $S_2(\Z)$ denote the group of cuspforms of weight~$2$ on~$\Gamma_0(N)$
with integral Fourier coefficients, and 
if $G$ is a subgroup of~$S_2(\Z)$, let
$G^\perp$ denote the subgroup of~$S_2(\Z)$ consisting
of cuspforms that are orthogonal to every~$g$ in~$G$ with respect
to the Petersson inner product.
The exponent of the quotient group
$$\frac{S_2(\Z)}{S_2(\Z)[I_f] + S_2(\Z)[I_f]^\perp}$$
 is called
the {\em congruence exponent} of~$A_f$ (really, that of~$f$), denoted~$\rAfe$
and its order is called the {\em congruence number}, denoted~$\rAf$
(see~\cite[\S3]{ars:moddeg}).
If $f$ has integer Fourier coefficients, so that $A_f$ is an elliptic curve,
then the quotient group above is a quotient of~$\Z$
(e.g., by~(\ref{eqn:congexp}) and~(\ref{eqn:congexp2}) in Section~\ref{sec:quotients}), 
so that $\rAf = \rAfe$,
and either of them
is the largest integer~$r$ such that
there exists a cuspform $g \in S_2(\Z)$ that is orthogonal to~$f$
under the Petersson inner product and whose $n$-th Fourier coefficient
is congruent modulo~$r$ to the $n$-th Fourier coefficient of~$f$
for all positive integers~$n$.
We say that a prime is a {\em congruence prime for~$A_f$} if it divides
the congruence number~$\rAf$.

Congruence primes have been studied by Doi, Hida, Ribet, 
Mazur, and others (see, e.g.,~\cite[\S1]{ribet:modp}), 
and played an important role in Wiles's work~\cite{wiles:fermat} 
on Fermat's last theorem.  Frey and Mai-Murty have 
observed that an
appropriate asymptotic bound on the modular degree is equivalent to
the $abc$-conjecture (see~\cite[p.544]{frey:boston}
and~\cite[p.180]{murty:congruence}). 
Thus congruence primes and the modular degree are quantities of
significant interest. 
Theorem~3.6 of~\cite{ars:moddeg} says that
the modular exonent~$\nAfe$ divides the congruence exponent~$\rAfe$ and
if $p$ is a prime such that $p^2 \ndiv N$,
then $\ord_p(\nAfe) = \ord_p(\rAfe)$. 

One might wonder if similar relations hold between the modular number~$\rAf$
and congruence number~$\nAf$ (as opposed to modular/congruence {\em exponents}). 
As mentioned earlier, if $A_f$ is an elliptic curve, then $\nAf = \nAfe^2$
and $\rAf = \rAfe$,
and so one sees that in this case, we have $\nAf \mid \rAfs$. 
So to start with, one might wonder whether $\nAf$ divides~$\rAfs$
even if $A_f$ is not an elliptic curve (i.e., has dimension more than one);
this question makes sense also because $\nAf$ is a perfect square, while
$\rAf$ need not be a perfect square.
It turns out that the answer to the question is no: as mentioned 
in~\cite[Remark~3.7]{ars:moddeg} we have

\begin{eg} \label{eg}
There is a newform
  of degree $24$ in $S_2(\Gamma_0(431))$ such that
 $$\nAf = (2^{11}\cdot 6947)^2 \,\,\ndiv\,\, r^2_{A_f} = (2^{10}\cdot
  6947)^2.$$
\end{eg}

We say that a maximal ideal~$\m$ of~$\T$ satisfies {\it multiplicity one}
if $J_0(N)[\m]$ is of dimension two over~$\T/\m$. 
The reason one calls this ``multiplicitly one'' is that if the canonical
two dimensional representation~$\rho_\m$ over~$\T/\m$ attached to~$\m$ 
(e.g., see~\cite[Prop.~5.1]{ribet:modreps}) is irreducible, then
$J_0(N)[\m]$ is a direct sum of copies of~$\rho_\m$
(e.g., see~\cite[Thm.~5.2]{ribet:modreps}), and 
a maximal ideal~$\m$ of~$\T$ satisfies multiplicity one
precisely if the multiplicity of~$\rho_\m$ in this decomposition is one.
Even if $\rho_\m$ is reducible, the definition of multiplicity one given
above is relevant (e.g., see~\cite[Cor.~16.3]{mazur:eisenstein}).
The notion of multiplicity one plays an important role in arithmetic
geometry, and in particular was used in the proof of Fermat's last theorem
in~\cite{wiles:fermat}. 

It was remarked in~\cite{ars:moddeg} that concerning Example~\ref{eg} above 
where $\nAf \ndiv \rAfs$, 
the level~$431$ is prime and 
by~\cite{kilford:non-gorenstein},  mod~$2$ multiplicity one fails 
for~$J_0(431)$.
In this article, we show that multiplicity one is the only obstruction
for the divisibility $\nAf \mid \rAfs$ to fail. In fact, we show something
stronger:

\begin{thm} \label{thm:main2}
Let $p$ be a prime such that 
every maximal ideal~$\m$ with residue characteristic~$p$
that contains~$I_f$ satisfies multiplicity one.
Then $\ord_p(\nAf)  = \ord_p(\rAfs)$.
\end{thm}

The theorem above follows from the more general Theorem~\ref{thm:main}
and the discussion preceding it in Section~\ref{sec:quotients}.
Example~\ref{eg} above shows that the multiplicity one hypothesis 
is cannot be removed.

The theorem above is the analog of
Proposition~5.9 of~\cite{ars:moddeg}, which says that under the hypotheses
of the theorem above, we have $\ord_p(\nAfe)  = \ord_p(\rAfe)$. 
If $A_f$ is an elliptic curve, then as remarked earlier,
$\nAf = \nAfe^2$ and $\rAf = \rAfe$, so our theorem adds nothing new.

In the context of Example~\ref{eg}, our theorem gives a new proof
that mod~$2$ multiplicity fails for~$J_0(431)$ (the original proof
being the one in~\cite{kilford:non-gorenstein}). 
Note that in~\cite{ars:moddeg}, the authors found examples 
of failure of multiplicity one using  Propostion~5.9 of loc. cit., which
implies that if the modular exponent does not equal 
the congruence exponent for 
some newform~$f$, then there is a maximal ideal of~$\T$ 
that not satisfy multiplicity one. However, we could not have
detected the failure of multiplicity one 
in Example~\ref{eg} by checking if the modular
{\em exponent} equals the congruence {\em exponent}, 
since the equality holds in the example for any newform~$f$
by~\cite[Thm.~3.6(b)]{ars:moddeg}, 
considering that the level is prime in the example. 
At the same time, consideration of the modular {\em number} and the
congruence {\em number} did dectect the failure of multiplicity one.
It would be interesting to do more calculations to see
when $\nAf \ndiv \rAfs$, as this may give new instances of failure
of multiplicity one.

We remark that our theorem gives information about the {\em order}
of a certain intersection of abelian subvarieties of~$J_0(N)$
(as opposed to just about primes that divide this order)
in terms of congruences between modular forms
(in fact, we give information in a more general setting in 
Section~\ref{sec:quotients}). We expect that 
the relation between a particular such intersection
and certain congruences will be useful in understanding the ``visible factor''
in~\cite{agmer},
and hope that such relations will be useful in other contexts as well.

It is known that multiplicity one holds in several situations.
We content ourselves by pointing out that by the main theorem
in Section~1.2 of~\cite{mazur-ribet}, a maximal ideal~$\m$
with residue characteristic~$p$ satisfies multiplicity one if
either $p \ndiv N$ or $p || N$ and $\rho_\m$ is not
modular of level~$N/p$.
We also have:

\begin{prop} \label{prop:ribet}
Let $p$ be an odd prime and $\m$ be a maximal ideal of~$\T$
with residue characteristic~$p$ such that $\rho_\m$ is irreducible.
Assume that either\\
(i) $p \ndiv N$ or \\
(ii) $p || N$ and $I_f \subseteq \m$ for some newform~$f$. \\
Then $\m$ satisfies multiplicity one.
\end{prop}
\begin{proof}
If $p \nmid N$, then the claim follows from
Theorem~5.2(b) of~\cite{ribet:modreps}, so let us assume that $p || N$.
Let $X_0(N)_{\Z_p}$ denote the minimal regular resolution of the compactified
coarse moduli scheme over~$\Z_p$ associated to~$\Gamma_0(N)$ as
in~\cite[\S~IV.3]{deligne-rapoport} and let 
$\Omega_{X_0(N)_{\Z_p}/\Z_p}$ 
denote the relative dualizing sheaf of~$X_0(N)_{\Z_p}$
over~$\Z_p$
(it is the sheaf of regular differentials as in~\cite[\S7]{mazur-ribet}).
We denote by $X_0(N)_{\Fp}$ the special fiber of $X_0(N)_{\Zlp}$ at the
prime~$p$ and by~$\Omega_{X_0(N)/\Fp}$ 
the relative dualizing sheaf of~$X_0(N)_{\Fp}$
over~$\Fp$.
It is shown in~\cite[\S5.2.2]{ars:moddeg} that under
the hypotheses above, 
${\rm dim}_{\T/\m} H^0(X_0(N)_{\Fp}, \Omega_{X_0(N)_{\Fp}/\Fp}) [\m] \leq 1$.
Let $J_{\Z_p}$ denote the N\'eron model of~$J_0(N)$ over~$\Z_p$
and let $J_{\Z_p}^0$ denote its identity component.
Then the natural morphism ${\rm Pic}^0_{X_0(N)/{\Zlp}} \ra J_{\Zlp}$
identifies ${\rm Pic}^0_{X_0(N)/\Zlp}$ with~$J_{\Zlp}^0$
(see, e.g., \cite[\S9.4--9.5]{neronmodels}).
Passing to tangent spaces along the identity section over~$\Zlp$,
we obtain an isomorphism
$H^1(X_0(N)_{\Zlp}, \OO_{X_0(N)_{\Zlp}}) \isom {\rm Tan}(J_{\Zlp})$.
Reducing both sides modulo~$p$ and applying Grothendieck
duality, we get 
${\rm Tan}(J_{\F_p}) \isom {\rm Hom}(H^0(X_0(N)_{\F_p}, \Omega_{X_0(N)/\Fp}), 
\Fp)$. Thus from the discussion above, we see that  
${\rm Tan}(J_{\F_p})/ \m {\rm Tan}(J_{\F_p})$ 
has dimension at most one over~$\T/\m$.
Since ${\rm Tan}(J_{\Zlp})$ is a faithful $\T \tensor \Z_p$-module,
we see that  
${\rm Tan}(J_{\F_p})/ \m {\rm Tan}(J_{\F_p})$ is non-trivial,
hence it is one dimensional
over~$\T/\m$. With this input, the proof of multiplicity one
in Theorem 2.1 of~\cite{wiles:fermat}, which is in the $\Gamma_1(N)$ context,
but is a formal argument involving 
abelian varieties (apart from the input above), 
carries over in the $\Gamma_0(N)$ context 
with the obvious modifications (in particular, replacing 
$X_1(N/p,p)_{\Z_p}$ in loc. cit. by~$X_0(N)_{\Z_p}$)
 to prove our claim (see p.~487-488 of loc. cit., as well 
as~\cite{tilouine:gorenstein}, where the input above
is the equation~(**) on p.~339).
\end{proof}

We remark that the condition that $p^2 \ndiv N$ in condition~(ii) 
of the proposition above cannot be removed, as follows from the counterexamples
in~\cite[\S2.2]{ars:moddeg}. From Theorem~\ref{thm:main2} 
and Proposition~\ref{prop:ribet}, we obtain:

\begin{cor}
Let $p$ be an odd prime.
Suppose that either \\
(i) $p \ndiv N$ or\\
(ii) $p || N$ and $\Afdual[\m]$ is irrreducible for
every maximal ideal~$\m$ of~$\T$ with residue characteristic~$p$.\\
Then $\ord_p(\nAf)  = \ord_p(\rAfs)$.
\end{cor}
\begin{proof}
The corollary is clear from
Theorem~\ref{thm:main2} and Proposition~\ref{prop:ribet}
in the case where $p \ndiv N$, so let
us assume that $p || N$. 
By Theorem~\ref{thm:main2} and Proposition~\ref{prop:ribet}, it suffices
to show that $\rho_\m$ is irreducible
for every maximal ideal $\m$ 
of~$\T$
with residue characteristic~$p$ such that $I_f \subseteq \m$.

Let $\m$ be such a maximal ideal. Then note that $\Afdual[\m]$ is non-trivial
since $\T/I_f$ acts faithfully on~$\Afdual$.
Let $D$ denote the direct sum of~$\Afdual[\m]$ and its Cartier dual. Let
$\ell$ be a prime that does not divide~$Np$ and let 
$\Frobl$ denote the Frobenius element of~${\rm Gal}(\Qbar/\Q)$ at~$\ell$.
As discussed in~\cite[p.~115]{mazur:eisenstein},
by the Eichler-Shimura relation, 
the characteristic polynomial of
$\Frobl$ acting on~$D$ is~$(X^2 -a_\ell X + \ell)^d =0$,
where $a_\ell$ is the image of~$T_\ell$ in~$\T/\m$ and
$d$ is the $\T/\m$-dimension of~$\Afdual[\m]$.
But this is also the characteristic polynomial of~$\Frobl$ acting
on the direct sum of $d$ copies of~$\rho_\m$. By the Chebotarev density theorem and the Brauer-Nesbitt theorem,
the semisimplification of~$D$ is $\rho_\m^d$. Thus
the semisimplification of~$\Afdual[\m]$ is a direct sum of certain number of
copies of~$\rho_\m$.
But $\Afdual[\m]$ is irreducible by hypothesis, so 
$\rho_\m = \Afdual[\m]$. Thus $\rho_\m$ is also irreducible, as was to be shown.
\end{proof}

Note that the hypothesis in the corollary above that
$\Afdual[\m]$ is irrreducible may be difficult to check, while it may be easier to check  
that $\Afdual[\m]$ does not have 
a nontrivial rational point (say in specific
instances, for example by checking that $p$ does not divide 
the order of the torsion subgroup of~$\Afdual(\Q)$); clearly the former 
hypothesis implies the latter, and for results in the other direction
(which is the relevant one in this context),
see~\cite{aw}. 

The corollary above
is the analog of Theorem~3.6(b) of~\cite{ars:moddeg}, which says
that $\ord_p(\nAfe) = \ord_p(\rAfe)$ provided $p^2 \ndiv N$,
in the setting of modular/congruence {\em numbers} as opposed
to modular/congruence {\em exponents} (although, note that 
we have an extra irreducibility hypothesis in our corollary).
We remark that the proofs of both results rely ultimately on ``multiplicity one
for differentials'' (as defined in~\cite[\S5.2]{ars:moddeg}).

If the level~$N$ is prime, then more can be said.
By Prop.~II.14.2 and Corollary~II.16.3 of~\cite{mazur:eisenstein}, 
every maximal ideal~$\m$ 
such 
that $\rho_\m$ is reducible also satisfies multiplicity one. 
Thus in view of 
Theorem~\ref{thm:main2} and Proposition~\ref{prop:ribet},
we obtain the following:
\begin{cor}
Suppose the level~$N$ is prime and let $p$ be an odd prime.
Then $\ord_p(\nAf)  = \ord_p(\rAfs)$.
\end{cor}

Also, much is known
in this situation if $\rho_\m$ is irreducible and $\m$ has residue characteristic~$2$ 
-- we refer to~\cite{kilford:non-gorenstein} and the references therein 
for details.
But note that by the examples in~\cite{kilford:non-gorenstein} or by
Example~\ref{eg} and Theorem~\ref{thm:main2}, multiplicity one need not
hold for a maximal ideal~$\m$ of residue characteristic~$2$ with
$\rho_\m$ irreducible even if the level~$N$ is prime.

The organization of the rest of this article is as follows.
In Section~\ref{sec:quotients},
we describe a more general setup, which includes
abelian subvarieties of~$J_1(N)$, and state a more 
general version of Theorem~\ref{thm:main2} (Theorem~\ref{thm:main} 
below). Section~\ref{sec:proof2} is devoted to the proof of 
the main result in Section~\ref{sec:quotients}.\\

\noindent {\bf Acknowldegements}: \\
The author is grateful to K.~Ribet for simplifying the proof
of Theorem~\ref{thm:main3} and for many suggestions that improved this paper.
We also thank J.~Tilouine for discussions regarding the proof of 
Proposition~\ref{prop:ribet}.\\

\section{A more general setup and some more results}
\label{sec:quotients}

For the benefit of the reader,
we first recall some of the discussion in~\cite[Section~3]{ars:moddeg}.
For $N\geq 4$, let~$\Gamma$ be either $\Gamma_0(N)$ or~$\Gamma_1(N)$.
Let~$X$ denote the modular curve over~$\Q$ associated
to~$\Gamma$, and let~$J$ be the Jacobian of~$X$.  
Let $J_f$ denote the standard abelian subvariety of~$J$
attached to~$f$ by Shimura~\cite[Thm.~7.14]{shimura:intro}.
Up to isogeny, $J$ is the product of factors $J_f^{e(f)}$ where $f$  
runs over the set of newforms of level dividing~$N$, taken up to Galois  
action, and $e(f)$ is the number of divisors of~$N/N(f)$, 
where $N(f)$ is the level of~$f$.  
Let $A$ be the sum of~$J_f^{e(f)}$ for some set of~$f$'s
(taken up to Galois action),
and let $B$ be the sum of all the other~$J_f^{e(f)}$'s.
Clearly $A + B = J$. The $J_f$'s are simple (over~$\Q$), 
hence $A \cap B$ is finite. 
By~\cite[Lemma~3.1]{ars:moddeg}, 
$\End(J)$ preserves $A$ and~$B$, 
where if $C$ is an abelian variety over~$\Q$,
by $\End(C)$ we mean the ring of endomorphisms
of~$C$ defined over~$\Q$.
If $f$ is a newform of weight~$2$ on~$\Gamma$ and $A_f$ is its
associated newform quotient, then $\Afdual$ and~$I_f J$
provide an example of $A$ and~$B$ respectively
as above, 
as shown in the discussion following Lemma~3.1 in~\cite{ars:moddeg}.

There is an alternate way to describe the~$A$ and~$B$ in the previous paragraph.
Since $\T \tensor \Q$ breaks up as a direct sum of algebras corresponding
to Galois orbits  of newforms of level dividing~$N$, the abelian subvariety~$A$ 
corresponds to an idempotent $e \in \T \tensor \Q$, and conversely, given 
an idempotent $e \in \T \tensor \Q$, the image of~$J$ under~$e$ (viewed 
as an element of~$\End(J) \tensor \Q$, which is to be multiplied by a large enough integer
to make the product an element of~$\End(J)$), is the corresponding~$A$ (and then
$B$ is the image of~$(1-e)$).

The {\em modular exponent}~$\nAe$ of~$A$ is 
defined as the exponent of~$A \cap B$
and the {\em modular number}~$\nA$ of~$A$ is its order
(see~\cite[\S3]{ars:moddeg}).
Note that the definition is symmetric with respect to~$A$ and~$B$.
If $f$ is a newform, then by
the modular exponent/number of~$A_f$, we mean that of~$A = \Afdual$, with
$B = I_f J$, which agrees with our earlier definition.

If~$R$ is a subring of~$\C$, 
let $S_2(R)=S_2(\Gamma;R)$ denote the subgroup of~$S_2(\Gamma; \C)$
consisting of cups forms whose Fourier expansions at the cusp~$\infty$
have coefficients in~$R$.  
Let $\T$ denote the Hecke algebra corresponding to the group~$\Gamma$.
There is a $\T$-equivariant bilinear pairing 
\begin{eqnarray} \label{eqn:pairing}
\T\times
S_2(\Z)\to\Z
\end{eqnarray}
given by $(t,g)\mapsto a_1(t(g))$, 
which is perfect 
(e.g., see \cite[Lemma~2.1]{abbull} or~\cite[Theorem~2.2]{ribet:modp}). 
Let $\T_A$ denote the image of~$\T$ in $\End(A)$, 
and let $\T_B$ be the image of $\T$ in $\End(B)$
(since $\T \subseteq \End(J)$, $\T$ preserves~$A$ and~$B$). 
Since $A + B = J$, the natural map $\T \ra \T_A \oplus \T_B$
is injective, and moreover, its cokernel is finite
(since $A \cap B$ is finite). 

Let $S_A = \Hom(\T_A,\Z)$ and~$S_B = \Hom(\T_B,\Z)$ be the
subgroups of~$S_2(\Z)$ obtained via the pairing in~(\ref{eqn:pairing}).
By~\cite[Lemma~3.3]{ars:moddeg}, we 
have an isomorphism
\begin{equation}\label{eqn:congexp}
   \frac{S_2(\Z)} { S_A + S_B} \isom \frac{\T_A \oplus \T_B}{\T}\ . 
\end{equation}
Also, we have an isomorphism
\begin{equation}\label{eqn:congexp2}
\frac{\T}{I_A + I_B} \stackrel{\simeq}{\longrightarrow}
\frac{\T_A \oplus \T_B}{\T}
\end{equation}
obtained by sending $t \in \T$ to~$(\pi_A(t),0) \in \T_A \oplus \T_B$,
where $\pi_A$ is the projection map $\T \ra \T_A$.
By definition~\cite{ars:moddeg}, 
the exponent of either of the isomorphic groups in~(\ref{eqn:congexp})
or~(\ref{eqn:congexp2})
is the {\em congruence exponent} $\rAe$ of~$A$ and the order of either
group is the {\em congruence number} $\rA$.
Note that this definition is also symmetric with respect to~$A$ and~$B$,
and again, the definition depends on both~$A$ and~$B$, unlike what 
the notation may suggest -- we have suppressed the dependence on~$B$
with the implicit understanding that $B$ has been chosen (given~$A$).
If $f$ is a newform, then by
the congruence exponent/number of~$A_f$, 
we mean that of~$A = \Afdual$, with
$B = I_f J$. In this situation,
$\T_A = \T/I_f$ and $S_A = S_2(\Z)[I_f]$. 
Also, $\Hom(\T_B,\Z)$ is the unique saturated
Hecke-stable complement of $S_2(\Z)[I_f]$ in $S_2(\Z)$, hence
must equal $S_2(\Z)[I_f]^{\perp}$. This shows that the new definition
of the congruence number/exponent generalizes our earlier 
definition for~$A_f$.

Let $I_A = \Ann_\T(A)$ and $I_B = \Ann_\T(B)$.
Theorem~3.6(a) of~\cite{ars:moddeg} says that 
the modular exponent~$\nAe$ divides the congruence exponent~$\rAe$, 
and Propostion~5.9 of loc. cit. says that if
$p$ is a prime such that all maximal ideals~$\m$ of~$\T$ containing
$I_A + I_B$ satisfy multiplicity one, then
$\ord_p(\rAe) = \ord_p(\nAe)$. 
Our main theorem deals with the case of modular/congruence numbers
as opposed to modular/congruence exponents. In view of the case
of newform quotients discussed in Section~\ref{sec:intro},
one would like to understand the relation between the modular
number~$\nA$ and the {\em square} of the congruence number~$\rA$.
As mentioned earlier,
it is not true that $\nA$ divides~$\rAs$ in general. At
the same time, we have:

\begin{thm} \label{thm:main}
Let $p$ be a prime  such that 
every maximal ideal of~$\T$ with residue characteristic~$p$
that contains~$I_A + I_B$ satisfies multiplicity one.
Then $\ord_p(\nA)  = \ord_p(\rAs)$.
\end{thm}

The theorem above follows immediately from:

\begin{thm}\label{thm:main3}
Let $\m$ be a  maximal ideal of~$\T$
that satisfies multiplicity one.
Then on tensoring with~$\T_\m$,
$A \cap B$ is free of rank two over
$\frac{\T}{I_A + I_B}$.
\end{thm}

We will prove this theorem in Section~\ref{sec:proof2}. 
Note that Theorem~\ref{thm:main}
is an analog of Propostion~5.9 of~\cite{ars:moddeg} in
the context of modular/congruence numbers as opposed to
modular/congruence exponents.
For results on
multiplicity one in the $\Gamma = \Gamma_1(N)$  context, 
see, e.g., \cite{tilouine:gorenstein} and the references therein.

Let $\pi_A:\T\to \T/I_A = \T_A$ and $\pi_B : \T \to \T/I_B = \T_B$ denote the natural
projection maps.
In this setup, in~\cite{ars:moddeg}, we defined
the {\em congruence ideal} as
$R=\pi_A(\ker(\pi_B)) \subset \T_A$, and 
the {\em intersection ideal} as
$  S = \Ann_{\T_A}(A \cap B)$. 
By~\cite[Lem.~5.2]{ars:moddeg}, we have $R \subseteq S$.
Moreover, $\pi_A$ induces a natural isomorphism
$$
\frac{\T}{I_A + I_B} \stackrel{\simeq}{\longrightarrow}
\frac{\T_A}{R}.$$
Thus from Thorem~\ref{thm:main3}, we obtain 

\begin{prop}
Let $\m$ be a  maximal ideal of~$\T$ 
that contains~$I_A$ and satisfies multiplicity one.
Then on tensoring with the completion of~$\m/I_A$, 
$R = S$. 
\end{prop}

The result above is not new: it follows from Proposition~5.6
and Lemma~5.8 of~\cite{ars:moddeg}. Our proof above is analogous
to the arguments in~\cite{ars:moddeg} (in this article, we use homology groups
in the proof of Thorem~\ref{thm:main3}, 
while the proof of Lemma~5.8 of~\cite{ars:moddeg} uses
Tate modules).

\section{Proof of Theorem~\ref{thm:main3}} \label{sec:proof2}
We have (e.g., by~\cite[Prop.~3.2]{agst:bsd}) the following natural isomorphism
of~$\T$-modules:
\begin{eqnarray} \label{eqn:2}
A\cap B \isom   \frac{H_1(J,\Z)}{H_1(A,\Z) + H_1(B,\Z)} \nonumber
\end{eqnarray}
Also, by~\cite[Lem~4.3]{agmer}, we have 
$H_1(A,\Z)  = H_1(J,\Z)[I_A]$ and $H_1(B,\Z) = H_1(J,\Z)[I_B]$.
Thus 
\begin{eqnarray} \label{eqn:3}
A\cap B  \isom  
\frac{H_1(J,\Z)}{H_1(J,\Z)[I_A] + H_1(J,\Z)[I_B]}. 
\end{eqnarray}

Recall that $\m$ is a maximal ideal of~$\T$ that satisfies multiplicity one.
Then 
$H_1(J,\Z) \tensor_\T \T_\m$ is free of rank two over~$\T_\m$,
by a standard argument due to Mazur (see~\cite[Lem.~II.15.1]{mazur:eisenstein}
or Corollary~(3) of Theorem~3.4 in~\cite{tilouine:gorenstein}).
Thus by~(\ref{eqn:3}), on tensoring with~$\T_\m$, $A\cap B$ is free
of rank two over $\frac{\T}{\T[I_A] + \T[I_B]}$.
The theorem now follows from the following lemma.

\begin{lem} \label{lem:end}
$\T[I_A] = I_B$ and $\T[I_B] = I_A$. 
\end{lem}

Before giving  the proof of this lemma, we need another lemma, which is
in a slightly more general setup.
In this paragraph, the symbol~$g$ stands
for a newform of some level~$N_g$ dividing~$N$. 
Let $S'_g$ denote the subspace of~$S_2(\Gamma_0(N),\C)$
spanned by the forms $g(dz)$ where $d$ ranges over the
divisors of~$N/N_g$. 
Let $[g]$ denote the Galois orbit of~$g$, 
and let $S_{[g]}$ denote the $\Q$-subspace of
forms in~$\oplus_{h \in [g]} S'_h$ with rational Fourier
coefficients. We have 
$S_2(\Gamma_0(N),\Q) = \oplus_{[g]} S_{[g]}$, where
the sum is over Galois conjugacy classes of newforms of
some level dividing~$N$.
Let $X$ be a subset of the set of Galois conjugacy classes of newforms of
some level dividing~$N$, and let 
$I = \Ann_\T (\oplus_{[g] \in X} S_{[g]})$. 
If $g$ is a newform of some level dividing~$N$,
then $S_{[g]}$ is preserved by~$\T$; let
$\T_{[g]}$ denote the image of~$\T$ acting on~$S_{[g]}$.
Then the natural map 
$$ \phi: \T \tensor \Q \ra \oplus_{[g]} \T_{[g]}$$
is an isomorphism of $\T \tensor \Q$ algebras, where $[g]$ ranges over all Galois conjugacy
classes of newforms of level dividing~$N$ (see, e.g., \cite[Thm.~3.5]{parent}). 
We have the decomposition
\begin{eqnarray}\label{eqn:decomp}
\oplus_{[g]} \T_{[g]} = \big(\oplus_{[g] \in X} \T_{[g]} \big)\  
\oplus\ \big(\oplus_{[g] \not\in X} \T_{[g]}\big) \ .
\end{eqnarray}
Let $\Ihat$ denote $\Ann_{\T}(I)$.

\begin{lem}\label{lem:decomp}
The image of $I \tensor \Q$
under $\phi$ is $\oplus_{[g] \not\in X} \T_{[g]}$,
and the image of~$\Ihat \tensor \Q$ is 
$\oplus_{[g] \in X} \T_{[g]}$. Thus
$\T \tensor \Q \isom I \tensor \Q \oplus \Ihat \tensor \Q$
as $\T \tensor \Q$-modules. 
\end{lem}
\begin{proof}
It is clear that the image of~$I \tensor \Q$ in~(\ref{eqn:decomp})
under $\phi$ is $\oplus_{[g] \not\in X} \T_{[g]}$.
As for the image of~$\Ihat \tensor \Q$, it clearly
contains $\oplus_{[g] \in X} \T_{[g]}$. Conversely, if 
$x \in \Ihat \tensor \Q$,  then it annihilates the element
$(0,1)$ in the decomposition of~(\ref{eqn:decomp}) (since
$(0,1)$ is
in the image of~$I \tensor \Q$ under~$\phi$), so the image of $x \cdot (0,1)$
in $\oplus_{[g] } \T_{[g]}$ must
be zero. Thus 
$x \in \oplus_{[g] \in X} \T_{[g]}$, which
finishes the proof the lemma.
\end{proof}

\begin{proof}[Proof of Lemma~\ref{lem:end}]
First, note that $\T[I_A] = \IAhat$.
By Lemma~\ref{lem:decomp}, taking $X$ to be the set of newforms corresponding
to~$A$, we see that the image of $I_A \tensor \Q$ under~$\phi$ in the 
decomposition~(\ref{eqn:decomp}) is $\oplus_{[g] \not\in X} \T_{[g]}$,
and by a similar argument, the image of $I_B \tensor \Q$
is $\oplus_{[g] \in X} \T_{[g]}$. Thus $I_B \subseteq \IAhat$ and
also, by Lemma~\ref{lem:decomp}, 
we see that $I_B \tensor \Q = \IAhat \tensor \Q$. But $I_B$ and~$\IAhat$ are
both saturated in~$\T$, and so it follows that $I_B = \IAhat$. This shows
that $I_B = \T[I_A]$. Swapping the roles of~$A$ and~$B$, we see that
$I_A = \T[I_B]$, which finishes the proof of the lemma.
\end{proof}

\bibliographystyle{amsalpha}         

\bibliography{biblio}

\end{document}